\documentclass[12pt,draft]{amsart}
\usepackage{amsmath,amsthm,latexsym,amscd,amsbsy,amssymb,graphics}
 \relax

\chardef\bslash=`\\ 

\makeatletter
\def\verbatim{\interlinepenalty\@M \@verbatim
  \leftskip\@totalleftmargin\advance\leftskip2pc
  \frenchspacing\@vobeyspaces \@xverbatim}
\makeatother
\makeatletter
  \def\dgt@k{\dg@DX=1 \dg@DY=-4 \dg@SIZE=1}
\makeatother

\theoremstyle{plain}

\newtheorem*{A}{Conjecture A}
\newtheorem*{B}{Conjecture B}

\theoremstyle{definition}

\newcounter{rmnum}

\numberwithin{equation}{section}

\begin{document}


\title[Notes on two conjectures in Extension Theory]
{Notes on two conjectures in Extension Theory}
\author{A.~Chigogidze}
\address{Department of Mathematics and Statistics,
University of Saskatche\-wan,
McLean Hall, 106 Wiggins Road, Saskatoon, SK, S7N 5E6,
Canada}
\email{chigogid@math.usask.ca}
\thanks{Author was partially supported by NSERC research grant.}



\begin{abstract}{It is noted that conjectures about the non-existence of universal compacta and compactifications of the given extension dimension  for non finitely dominated complexes are not valid for all CW complexes of the form $L \vee S^{2}$, where $L$ is of finite type and has a finite fundamental group, but is not finitely dominated.}
\end{abstract}

\maketitle \markboth{A.~Chigogidze}{On two conjectures of J.~Dydak}

\section{Introduction}\label{S:intro}

The following conjecture has been stated in \cite[Conjecture 2.8]{dydak1}
and later restated in \cite[Conjecture 1.23]{DD}:

\begin{A}\label{A}
Suppose $K$ is a countable CW complex. If the class of compacta
$\{ X \colon K \in AE(X)\}$ has a universal space, then $K$ is
homotopy dominated by a finite CW complex.
\end{A}

Closely related to this is the following hypothesis \cite[Conjecture 1.20]{DD}:

\begin{B}\label{B}
Let $K$ be a countable CW complex. Any separable metrizable space $X$
with $K \in AE(X)$ admits a metrizable compactification $\widetilde{X}$
with $K \in AE(\widetilde{X})$ if and only if $K$ is homotopy dominated
by a finite CW complex.
\end{B}

Below we show that both conjectures fail (see \cite[Lemma 1.11]{lrs} where Conjecture A has been formally disproved) for a large class of complexes. This becomes clear by observing that
the assertions of these conjectures are valid for arbitrary CW complexes of the form
$K = L\vee S^{2}$, where $L$ is a CW complex of finite type which is not finitely dominated and which has a finite
fundamental group (consider, for instance, $L = {\mathbb R}P^{\infty}$). Note that $K$ also is not finitely dominated and has a finite fundamental group. 

Author is grateful to J.Dydak for helpful e-mail exchange.

\section{Disproving Conjectures A and B}\label{S:disp}

For spaces $X$ and $L$ the notation $L \in {\rm AE}(X)$ means that every continuous map $f \colon A \to L$, defined on a closed subspace $A$ of $X$, admits a continuous extension $\widetilde{f} \colon X \to L$.  A.~N.~Dranishnikov introduced \cite{dra} the following partial order for CW complexes. We say that $L \leq K$ if for each space $X$ the condition $L \in {\rm AE}(X)$ implies the condition $K \in {\rm AE}(X)$. Equivalence classes of complexes with respect to this relation are called {\em extension types}. The above defined relation $\leq$ creates a partial order in the class of extension types of complexes. This partial order is still denoted by $\leq$ and the extension type with representative $K$ is denoted by $[K]$ (\cite{chi3}, where the reader can find a general overview of extension theory, discusses further properties of extension types in terms of this partial order).

For the reader's convenience we state here two statements which are needed below. Recall that a complex has a finite type if it contains finite number of cells in every dimension. $\beta X$ denotes the Stone-\v{C}ech compactification of a space $X$.

{\rm {\bf Proposition 1} \cite[Proposition 2.3]{chi2}. {\em Let $L$ be a CW complex of finite type with a finite fundamental group. Then the following conditions are equivalent for any normal space $X$ and any integer $n \geq 2$:
\begin{itemize}
\item[(a)]
$L \vee S^{n} \in {\rm AE}(X)$.
\item[(b)]
$L \vee S^{n} \in {\rm AE}(\beta X)$.
\end{itemize}}

{\rm {\bf Proposition 2} \cite[Theorem 2.1]{chi2}. {\em Let $L$ be a countable CW complex. Then the following conditions are equivalent for any normal space $X$:
\begin{itemize}
\item[(a)]
$L \in {\rm AE}(\beta X)$ whenever $L \in {\rm AE}(X)$.
\item[(b)]
There exists an $L$-invertible map $f_{L} \colon X_{L} \to Q$ of a metrizable compactum $X_{L}$, with $L \in {\rm AE}(X_{L})$, onto the Hilbert cube $Q$.
\end{itemize}}

\subsection{Conjecture A}
In order to prove our claim related to Conjecture A, take a complex $K = L \vee S^{2}$ of the type indicated in the Introduction and note that, by Proposition 1, $K \in {\rm AE}(\beta Y)$ for any normal space $Y$ with $K \in {\rm AE}(Y)$. Then, by Proposition 2, there exists a $K$-invertible map $f_{K} \colon Y_{K} \to Q$ of a metrizable compactum $Y_{K}$, with $K \in {\rm AE}(Y_{K})$, onto the Hilbert cube $Q$. This implies that any embedding $Z \hookrightarrow Q$ of a metrizable compactum  with $Z \in {\rm AE}(X)$ can be lifted to an embedding $Z \hookrightarrow X_{K}$. Consequently $Y_{K}$ is an universal metrizable compactum for the class indicated in Conjecture A.

\subsection{Conjecture B}
Next consider a separable metrizable space $X$ such that $K \in {\rm AE}(X)$. By \cite{ols}, we may assume that $X$ is completely metrizable. By Proposition 1, $K \in {\rm AE}(\beta X)$. By \cite[Theorem 4.4]{chi1}, $\beta X$ can be represented as the limit of a Polish spectrum ${\mathcal S} = \{ X_{\alpha}, p_{\alpha}^{\beta}, A\}$ consisting of Polish spaces $X_{\alpha}$ such that $K \in {\rm AE}(X_{\alpha})$, $\alpha \in A$. Since $\beta X$ is compact, we may assume that all $X_{\alpha}$'s are metrizable compact spaces. We will show that there exists an index $\widetilde{\beta} \in A$ such that the restriction $p_{\widetilde{\beta}}|X \colon X \to p_{\widetilde{\beta}}(X)$ of the limit projection $p_{\widetilde{\beta}} \colon \beta X \to X_{\widetilde{\beta}}$ of the spectrum ${\mathcal S}$ is a homeomorphism. This would complete the proof since the metrizable compactum $\widetilde{X} = X_{\widetilde{\beta}}$ would then serve as the required compactification of the space $X$.

Since $X$ is completely metrizable there exist functionally open subsets $G_{k}$ in $\beta X$ such that $X = \cap_{k=1}^{\infty}G_{k}$. Since the spectrum ${\mathcal S}$ is factorizing $\sigma$-spectrum (see \cite[Sections 1.3.1 and 1.3.2]{chibook}) we can conclude that for each $k$ there exists an index $\alpha_{k} \in A$ such that $G_{k} = p_{\alpha_{k}}^{-1}\left( p_{\alpha_{k}}(G_{k})\right)$. According to \cite[Corollary 1.1.28]{chibook} there exists an index $\widetilde{\alpha} \in A$ such that $\widetilde{\alpha} \geq \alpha_{k}$ for each $k$. Then $X = p_{\widetilde{\alpha}}^{-1}\left( p_{\widetilde{\alpha}}(X)\right)$. Moreover, it is easy to see that $X = p_{\alpha}^{-1}\left( p_{\alpha}(X)\right)$ for each $\alpha \geq \widetilde{\alpha}$. It then follows that for each such $\alpha$ the restriction $p_{\alpha}|X \colon X \to p_{\alpha}(X)$ is a perfect surjection.

Finally we show that there exists an index $\widetilde{\beta} \geq \widetilde{\alpha}$ such that every its fiber $\left( p_{\widetilde{\beta}}|X\right)^{-1}(y)$, $y \in p_{\widetilde{\beta}}(X)$, consists of precisely one point. Choose a countable open base $\{ U\}_{k=1}^{\infty}$ of $X$ (such a base exists since $X$ is a separable metrizable space). Let $V_{k}$ be a functionally open subset of the Stone-\v{C}ech compactification $\beta X$ such that $V_{k} \cap X = U_{k}$. Since the spectrum ${\mathcal S}$ is factorizing, there exists an index $\beta_{k} \geq \widetilde{\alpha}$ such that $V_{k} = p_{\beta_{k}}^{-1}\left( p_{\beta_{k}}(V_{k})\right)$. By \cite[Corollary 1.1.28]{chibook}, there exists an index $\widetilde{\beta} \in A$ such that $\widetilde{\beta} \geq \beta_{k}$ for each $k$. Let us show that this index has the required property. Assume the contrary, i.e. suppose that the exists a point $y \in p_{\widetilde{\beta}}(X)$ such that the fiber $\left( p_{\widetilde{\beta}}|X\right)^{-1}(y) = p_{\widetilde{\beta}}^{-1}(y) \cap X$ contains at least two distinct points. Denote them by $x_{1}$ and $x_{2}$ respectively. Since $\{ U_{k}\}_{k=1}^{\infty}$ is an open base of $X$ we can find an index $n$ such that $x_{1} \in U_{n}$ and $x_{2} \not\in U_{k}$. Note that 
\begin{multline*} U_{k} = V_{k} \cap X = p_{\beta_{k}}^{-1}\left( p_{\beta_{k}}(U_{k})\right) \cap X = p_{\widetilde{\beta}}^{-1}\left( \left( p_{\beta_{k}}^{\widetilde{\beta}}\right)^{-1}\left( p_{\beta_{k}}^{\widetilde{\beta}}( p_{\widetilde{\beta}}(U_{k})\right)\right) \cap X = \\
p_{\widetilde{\beta}}^{-1}\left( p_{\widetilde{\beta}}(U_{k})\right) \cap X .
\end{multline*}

\noindent Then, by the choice of the set $U_{k}$ we obtain two contradictory conclusions: 
\[ y = p_{\widetilde{\beta}}(x_{1}) \in p_{\widetilde{\beta}}\left( p_{\widetilde{\beta}}^{-1}\left( p_{\widetilde{\beta}}(U_{k})\right) \cap X\right) = p_{\widetilde{\beta}}(U_{k}) \cap p_{\widetilde{\beta}}(X)\]

\noindent and

\[ y = p_{\widetilde{\beta}}(x_{2}) \not\in p_{\widetilde{\beta}}\left( p_{\widetilde{\beta}}^{-1}\left( p_{\widetilde{\beta}}(U_{k})\right) \cap X\right) = p_{\widetilde{\beta}}(U_{k}) \cap p_{\widetilde{\beta}}(X) .\]

\noindent This shows, as required, that the restriction $p_{\widetilde{\beta}}|X \colon X \to p_{\widetilde{\beta}}$ is a homeomorphism.

\subsection {Concluding remarks}

CW complexes we used to disprove Conjectures A and B are of the form $K = L \vee S^{n}$, $n \geq 2$, where $L$ has a finite type (and a finite fundamental group). Of course, such complexes need not be finitely dominated. However, $[K] = [L \vee S^{n}] \leq [S^{n}]$ and consequently, according to \cite[Example 2.4(iv)]{chi3}, $[K] = [K^{(n)} \vee S^{n}] = [\left( L \vee S^{n}\right)^{(n)} \vee S^{n}] = [ L^{(n)} \vee S^{n}\vee S^{n}] = [L^{(n)}\vee S^{n}]$. This shows that the extension type $[K]$ contains a finite complex (note that $L^{(n)}$ is finite). This leads to a more plausible versions of the above conjectures. Namely, is it true that under the same assumptions it follows that the extension type $[K]$ of the complex contains a finitely dominated complex?


\begin{thebibliography}{99}
\bibitem{chi3}
A.~Chigogidze, {\em Infinite dimensional topology and shape theory}, in:
Handbook of Geometric Topology, ed. by R. Daverman and R. Sher, North Holland, Amsterdam, 2000, 307--371.


\bibitem{chi2}
A.~Chigogidze, {\em Compactifications and universal spaces in extension
theory}, Proc. Amer. Math. Soc. {\bf 128} (1999), 2187--2190.


\bibitem{chi1}
A.~Chigogidze, {\em Cohomological dimension of Tychonov spaces}, Topology Appl. {\bf 79} (1997), 197--228.

\bibitem{chibook}
A.~Chigogidze, {\em Inverse Spectra}, North Holland, Amsterdam, 1996.

\bibitem{dra}
A.~N.~Dranishnikov, {\em The Eilenberg-Borsuk theorem for maps in an arbitrary complex}, Mat. Sbornik {\bf 185} (1994), 81--90.


\bibitem{DD}
A.~Dranishnikov, J.~Dydak, {\em Extension theory of separable metric spaces
with applications to dimension theory}, Trans. Amer. Math. Soc. {\bf 353}
(2001), 133--156.


\bibitem{dydak1}
J.~Dydak, {\em Cohomological dimension of metrizable spaces. II}, Trans.
Amer. Math. Soc. {\bf 348} (1996), 1647 -- 1661.

\bibitem{lrs}
M.~Levin, L.~R.~Rubin, P.~J.~Schapiro, {\em The Marde\v{s}i\'{c} factorization theorem for extension theory and C-spaces}, Proc. Amer. Math. Soc. {\bf 128} (2000), 3099--3106.

\bibitem{ols}
W.~Olszewski, {\em Completion theorem for cohomological dimension}, Proc. Amer. Math. Soc. {\bf 123} (1995), 2261--2264.
\end{thebibliography}
\end{document}